\documentclass[12pt,reqno]{amsart}

\usepackage{anysize}

\marginsize{1in}{1in}{1in}{1in}

\usepackage{latexsym}
\usepackage{anyfontsize} 
\usepackage{tikz} 
\usetikzlibrary{cd} 
\usepackage{color}
\usepackage{hyperref}
\usepackage{url}
\usepackage{breakurl}
\usepackage{spverbatim}
\usepackage{ragged2e}
\newcommand{\bburl}[1]{\textcolor{blue}{\url{#1}}}

\makeatletter
\newcommand{\monthyear}[1]{%
  \def\@monthyear{\uppercase{#1}}}
\newcommand{\volnumber}[1]{%
  \def\@volnumber{\uppercase{#1}}}
\makeatother

\newcommand{\erf}{{\rm erf}}

\theoremstyle{plain}
\numberwithin{equation}{section} 
\newtheorem{thm}{Theorem}[section] 

\newtheorem{rek}[thm]{Remark}

\numberwithin{table}{section} 
\numberwithin{figure}{section}

\begin{document}

\setcounter{page}{1}

\title{An ODE to Changing Variables: $v = y/x$ Is For Victory}

\author[J. Kraisler]{Joseph Kraisler}
\email{\textcolor{blue}{\href{mailto:jk34@williams.edu}{jk34@williams.edu}}}
\address{Department of Mathematics, Williams College, Williamstown, MA 01267, USA}

\author[K. S. Miller]{Kayla S. Miller}
\email{\textcolor{blue}{\href{mailto:ksmill5433@gmail.com}{ksmill5433@gmail.com}}}
\address{Department of Mathematics, Williams College, Williamstown, MA 01267, USA}

\author[S. J. Miller]{Steven J. Miller}
\email{\textcolor{blue}{\href{mailto:sjm1@williams.edu}{sjm1@williams.edu},
\href{mailto:Steven.Miller.MC.96@aya.yale.edu}{Steven.Miller.MC.96@aya.yale.edu}}}
\address{Department of Mathematics, Williams College, Williamstown, MA 01267, USA}

\maketitle

\bigskip

\begin{abstract} Most first courses in differential equations introduce students to very special differential equations which can be solved in closed form from a small set of techniques (such as separable, integrable). Frequently other equations can be converted to one of these forms through a change of variable; however, these transformations are usually just given without any explanation of \emph{how} people found them, or \emph{why} problems of these forms have nice solutions, but not related ones. We consider the first order differential equation $dy/dx = f(x,y)$ where $f(x,y)$ doesn't factor as $A(x)B(y)$, and assume $f(x,y) = g(y/x^r)$ for some nice function $g$ and some real number $r$. This can be reduced to a separable one only when $r=1$. We then discuss generalizations and what conditions of $f$ lead to a change of variables yielding a separable equation. Our goals are twofold. First, to encourage students and instructors to be active, questioning the methods and results, seeing how far they can be extended and why they work. Second, we discuss how to appropriately use AI to take ideas and special cases as prompts and quickly create extensions and editable notes that can be shared with students, freeing up the professors (and students) from straightforward calculations \emph{without} sacrificing learning and understanding.
\end{abstract}


\section{Introduction}

Students frequently meet differential equations in their first year of calculus, studying simple problems (such as bacteria growth) with easily computed closed form solutions. Not surprisingly the situation changes quickly when they take a class devoted to the subject. Similar to the arsenal of integration techniques learned in calculus, they now see that equations of certain forms can be ``solved'' (see Appendix \ref{sec:appendixwhatsasoln} for a brief discussion on when we consider something ``solved'').

We study first order differential equations of the form $dy/dx = f(x,y)$.
This of course leads to the following, natural question: how to convert the equation we need to study to one of the forms we know how to study? For example, an equation is \emph{separable} if $f(x,y) = A(x) B(y)$, as in that case we have $$\int \frac{dy}{B(y)} \ = \int \ A(x) dx.$$

While we may not have nice closed forms for the resulting integrals, there are at least simple numerical techniques or frequently series expansions that are available, and we'll thus consider any problem that can be reduced to a separable one a success. Even determining whether or not a given first order differential equation is separable can be a difficult task for undergraduate students, and systematic techniques for answering this question has been explored in previous expository articles \cite{Ci,Sc}.\footnote{In \cite{Sc} it's proved for an open convex $D$ with $f: D \subset \mathbb{R}^2 \to \mathbb{R}$ that $f(x,y) = \phi(x)\psi(y)$ for all $(x,y) \in D$ if and only if
$f(a,b)f(x,y) = f(x,b)f(a,y)$ for all $(a,b), (x,y)$. In \cite{Ci}, for the same $D$ if $f$, $f_x$, $f_y$ and $f_{xy}$ are continuous and $f(x,y)$ is never $0$ in $D$ and
$f(x,y)f_{xy}(x,y) = f_x(x,y)f_y(x,y)$ for all $(x,y) \in D$ then there are continuously differentiable functions $\phi(x)$ and $\psi(y)$ such that
$f(x,y) = \phi(x)\psi(y)$ for all $(x,y) \in D$.}
Our paper tackles a related question: When can an ordinary differential equation be made seprable through a change of variables?

Our first result concerns $f(x,y)$ of the form $g(y/x^r)$. Note assuming we're of this form excludes two extreme cases, namely $f(x,y)$ is a function only of $x$ or only of $y$, but of course in those cases it's already separable and there's nothing to do! If we had $f(x,y)$ is a function only of $y^a/x^b$ with $a, b$ non-zero, setting $r = b/a$ we see $y^a/x^b = (y/b^r)^a$, and thus $f$ depends only on $y/x^r$. There is thus no loss in looking at functions of $y/x^r$ over the more general case $y^a/x^b$.

\begin{rek} We could of course have chosen to write $v = x^s y$ and would then find the only tractable case is $s = -1$. We prefer to write the change of variable as $v = y/x^r$ as this aligns with the formulation seen in books, but of course it doesn't matter. The key idea is to have a free parameter, and explore what choices lead to good outcomes. \end{rek}

\begin{thm}\label{thm:mainxyr} Assume $dy/dx = f(x,y)$ with $f(x,y) = g(y/x^r)$ for some $r \neq 0$, and $f(x,y)$ does not factor as $A(x) B(y)$ (in that case the equation is already separable). Then the differential equation is separable for all $g$ only if $r=1$. Explicitly, the change of variable $v = y/x$ works when $r=1$, and there's no similar choice that works if $r=-1$ (i.e., the case $f(x,y)$ is a function of just the product $xy$).
\end{thm}

Building on this success, it's natural to ask what happens for more complicated relations. Is there a nice change of variable if perhaps $f(x,y)$ depends only on $y/e^x$? Or perhaps $\sin(y)/\cos(x)$? A more careful argument yields the following.

\begin{thm}\label{thm:maingeneralcase}
Suppose that the right hand side of the differential equation $dy/dx = f(x,y)$ can be rewritten $f(x,y) = g(h(x,y))$ for $g$ which is continuously differentiable on an interval $J$ and $h$ which is twice continuously differentiable. Then the substitution $v =h(x,y)$ leads to a separable differential equation for every $h$ if and only if $h$ has one of the two following forms:
\begin{align}
    h(x,y) = H(kx+y)\quad \text{ or }\quad g(x,y) = H\left(\frac{y-y_0}{x-x_0}\right),
\end{align}
with constants $k,x_0,y_0$ where $H$ is twice continuously differentiable and $H'\neq 0$ on $J$. We note that the particular form of the differential equation obtained by this substitution is heavily dependent on the function $h$, and that additional substitutions may work for some, but not all, functions $g$.
\end{thm}

We first prove Theorem \ref{thm:mainxyr} in \S\ref{sec:proofmainyxr}, and then discuss how to guide students in proving its extension in \S\ref{sec:proofmaingeneralcase}.

The genesis of this note is a suggested problem by the first named author to the second named author, a student in his Math 309 Differential Equations and Non-linear Dynamics class at Williams in Fall '26: solve \begin{equation} \frac{dy}{dx} \ = \ \frac{x^2 + xy + y^2}{x^2}. \end{equation} Times have changed greatly since the third named author took such a class at Yale in the mid-90's. Then, a huge part of the class emphasized techniques, including how to convert a problem to a form in a reference book. As computer programs and online resources have advanced magnitudes since then, and since many real world equations don't have nice closed form solutions, classes emphasize specific techniques less and spend more on applications (such as non-linear dynamics). Giving additional problems such as this, which are not graded, is a great way to alert students to the methods and techniques that exist while still allowing time to get to the interesting applications. For this problem, she shared the problem with the third named author. He remarked that when he was a student, the class used the fifth edition of \cite{BoDi} (\emph{his} father used the first edition, so keeping it in the family), and there was a section that if $f(x,y)$ can be written as a function of just $v = y/x$ then one can transform to a separable equation. This prompted the following excellent questions from the second to the third author: When else can you change variables and get an equation you can solve? Can you do it for $v = xy$? The two of them explored the case of $f(x,y)$ a function of $y/x^r$, and she reported back to her professor, who then worked with AI to generalize the argument to find general conditions for when $f(x,y)$ leads to a separable equation. We hope this note helps other instructors and students in seeing the value of asking questions and actively engaging in the material.

It also leads us to the second reason for writing this argument. In addition to highlighting the importance of asking good questions and exploring the consequences, and not just passively accepting the math that comes one's way, this provides an excellent opportunity to discuss how to use AI well in a classroom. Extending the argument is straightforward and almost surely in the literature, though not necessarily written in a friendly way for beginning students. By feeding good prompts one can quickly generate high quality notes, complete with good citations to the literature. This frees up time for other activities, but not at the cost of understanding.

\section{The Standard Case: Why only $y/x$?}\label{sec:proofmainyxr}

\subsection{A Standard Example}

Let's explore the suggested problem, where $f(x,y) = (x^2 + xy + y^2)/x^2$. Thus we wish to solve \begin{equation}\label{eq:suggestedprobleminsubsec} \frac{dy}{dx} \ = \ \frac{x^2 + xy + y^2}{x^2}. \end{equation} Recalling the sage advice of the great New England transcendentalist\footnote{A different type of transcendentalism then one normally encounters in mathematics, but his advice is a universal truth.} Henry David Thoreau, our first thought whenever we see something should be: \emph{Simplify, simplify}.\footnote{For example, if we need to take the derivative of $$(1701 x^{24} + 601 x) (5 x^3 + 2 x + 1)$$ it's natural to use the product rule, which yields $$(1701 \cdot 24 x^{23} + 601) (5x^3 + 2x + 1) + (1701 x^{24} + 601 x) (5 \cdot 3 x^2 + 2). $$ However, rather than use the product rule if we first expand we get $$(1701 x^{24} + 601 x) (5 x^3 + 2 x + 1) \ = \ 8505 x^{27}+3402 x^{25}+1701 x^{24}+3005 x^4+1202 x^2+601 x, $$ and we can now easily take the derivative term by term using the sum and constant rules and thus avoid having to use the more involved product rule.}

How does that help us here? Simplifying the right hand side of \eqref{eq:suggestedprobleminsubsec} by doing the division gives $$\frac{x^2 + xy + y^2}{x^2} \ = \ 1 + \frac{y}{x} + \left(\frac{y}{x}\right)^2.$$ Thus $f(x,y)$ depends on $x$ and $y$ very weakly -- it only depends on their \emph{ratio}. This \emph{suggests} that it might be worthwhile changing variables and letting $v = y/x$, as now the right hand side is only explicitly a function of one variable, $v$ (as $y = y(x)$ of course $v = y/x$ means $v$ depends on $x$; this dependence is often suppressed in expositions, and textbooks frequently talk about the variables $x$ and $y$ or after the transformation $x$ and $v$). We now have \begin{equation} \frac{dy}{dx} \ = \ 1 + v + v^2. \end{equation} To make further progress, we need to replace $dy/dx$ with an expression that's a function of just $x, v$ and $dv/dx$. As $v = y/x$ (or $v(x) = y(x)/x$), the quotient rule gives \begin{equation} \frac{dv}{dx} \ = \ \frac{\frac{dy}{dx}x - y}{x^2}. \end{equation} Solving for $dy/dx$ gives \begin{equation} \frac{dy}{dx} \ = \ \frac{x^2 \frac{dv}{dx} + y}{x}, \end{equation} and then recalling that $v = y/x$ yields \begin{equation} \frac{dy}{dx} \ = \ x \frac{dv}{dx} + v. \end{equation}

Substituting for $dy/dx$ into our problem gives \begin{equation} x \frac{dv}{dx} + v \ = \ 1 + v + v^2 \ \ \ {\rm or} \ \ \ x \frac{dv}{dx} \ = \ 1 + v^2. \end{equation} The above is separable! After dividing we get \begin{equation} \frac{dv}{1 + v^2} \ = \ \frac{dx}{x}. \end{equation} Not only is our equation separable, but both integrals are relatively nice. The right hand side integrates to the well-known $\ln x$, while the left hand integrand is the less remembered (but quite important) anti-derivative of $\arctan$. Thus the solution is \begin{equation} \arctan(v) \ = \ \ln x + C \ \ \ {\rm or,\ since\ } v \ = \ y/x, \ \ \ \arctan(y/x) \ = \ \ln x + C, \end{equation} where $C$ is determined by the initial conditions.

Our problem had a nice cancellation: when replacing $dy/dx$ we had $x dv/dx + v$, and the $v$ cancelled with a $v$ on the right hand side. If instead of $1 + v + v^2$ on the right we had \emph{any} function of $v$, we would still have obtained a separable equation; all that matters is the right hand side depended only on $v$, and thus it could easily absorb the $v$ term from the left. Mark Twain has sage advice for situations such as these. \begin{quote} \emph{Good judgement is the result of experience and experience the result of bad judgement. We should be careful to get out of an experience only the wisdom that is in it and stop there lest we be like the cat that sits down on a hot stove lid.} \end{quote}

\subsection{General Power Relation}

We now prove Theorem \ref{thm:mainxyr}. Thus $dy/dx = f(x,r)$ and there's a non-zero $r$ and a nice function $g$ such that $f(x,y) = g(v)$ for $v = y/x^r$. Obviously if $f(x,y) = A(x) B(y)$ then the equation is already separable, so there's no need to further explore. The consequences of that observation for us is that if $g(v) = c v$ for some constant $c$ then our initial equation is already separable, as we then have $dy/dx = (c/x^r) \cdot y$, leading to \begin{equation} \frac{dy}{y} \ = \ \frac{c dx}{x^r}. \end{equation}

As before, the first step is to replace $dy/dx$ with an expression involving $dv/dx$. Using the quotient rule to differentiate $v = y/x^r$ (which we can view as $v(x) = y(x)/x^r$) we have \begin{equation} \frac{dv}{dx} \ = \ \frac{\frac{dy}{dx} x^r - y r x^{r-1}}{x^{2r}} \ = \  \frac{1}{x^r} \frac{dy}{dx} - \frac{r}{x^r} \frac{y}{x^r} x^{r-1}. \end{equation} We can simplify the above. On the right hand side we can replace $y/x^r$ with $v$, and then multiplying both sides by $x^r$ yields \begin{equation} x^r \frac{dv}{dx} \ = \ \frac{dy}{dx} - r v x^{r-1}, \ \ \ {\rm or} \ \ \ \frac{dy}{dx} \ = \ x^r \frac{dv}{dx} + rvx^{r-1}. \end{equation}

We substitute this into $dy/dx = f(x,y) = g(v)$ and find \begin{equation} x^r \frac{dv}{dx} + rvx^{r-1} \ = \ g(v), \ \ \ {\rm or} \ \ \ x^r \frac{dv}{dx} \ = \ g(v) - r v x^{r-1}. \end{equation}
The above is separable if and only if the right hand side factors as a function of $x$ times a function of $v$, say $\Phi(x) \Psi(v)$, as in that case we obtain \begin{equation} x^r \frac{dv}{dx} \ = \ \Phi(x) \Psi(v) \ \ \ {\rm or} \ \ \ \frac{dv}{\Psi(v)} \ = \ \frac{\Phi(x) dx}{x^r}.\end{equation}

So we need to find functions $\Phi$ and $\Psi$ such that \begin{equation} g(v) - r v x^{r-1} \ = \ \Phi(x) \Psi(v). \end{equation} If $r=1$ we can take $\Phi(x) = 1$ and $\Psi(v) = g(v) - v$. The other possibility is if $g(v) = c v$ for some constant $c$, as then \begin{equation} g(v) - r v x^{r-1} \ = \ c v - r v x^{r-1} \ = \ \left(c - r x^{r-1}\right) v. \end{equation}  However, we've already excluded that case as we're assuming $f(x,y)$ does not factor as $A(x) B(y)$, which we showed above would be the case if $g(v) = c v$.

Thus, if $f$ doesn't factor, only when $r=1$ will $g(v) - r v x^{r-1}$ have no explicit $x$-dependence and factor as needed. In that case, and that case only, we can then re-arrange and find \begin{equation} \frac{\ dv}{g(v) - v} \ = \ \frac{dx}{x}. \end{equation} Letting $\mathcal{G}$ be any anti-derivative of $1/(g(v) - v)$ we find \begin{equation} \mathcal{G}(v) \ = \ \ln x + C, \end{equation} solving the problem.

\section{The General Case: When $f$ Leads to Separable}\label{sec:proofmaingeneralcase}

As one of the main reasons for writing this article is to encourage students and teachers to actively engage with mathematics, we leave the proof of the general case to the reader. We've provided a hint by writing down the answer, but of course feel free to ignore and try to discover on your own.

We end with a few suggestions. First, it greatly helps to explore first to build intuition and understand the problem. We started with a specific problem, and then extended from a specific $f(x,y)$ that's a function of $y/x$ to looking at situations that are functions of $y/x^r$. Building on this, try a more general relation: $f(x,y) = f(h(x,y))$, and try to find conditions on $h$ that lead to a seprable equation.

Such investigations provide an outstanding opportunity to use AI in a good way. We are not outsourcing understanding; we've done some cases actively first, and can now prompt AI to try to push these further. Thus one approach is to write up some cases, and feed those arguments into AI and ask it to try to extend. While for a simple problem like this there's almost surely no need to give AI a good starting point by including our initial work, this is a good habit to get in. Namely, by doing some work first we can use AI to help us see just how far an idea can be extended and pushed.

\appendix
\section{What's a Solution?}\label{sec:appendixwhatsasoln}

We have to be careful about what we mean when we say we can solve something. Most students leave Calculus II believing they can integrate. Sadly most functions don't have nice closed-form anti-derivatives. However, if we allow ourselves to ``name'' the anti-derivative, then yes we can say many of these are solved! We give an example from probability, which is based on the fact that while there's no nice anti-derivative of $e^{-t^2}$, there is a rapidly converging series expansion.

We say $f$ is a probability distribution if it's non-negative and  integrates to 1; thus if $f$ is the density for a random variable $X$ then the probability $X$ takes on a value in $[a, b]$ is just $\int_a^b f(x) dx$.
We frequently write $F(x)$ for the associated cumulative distribution function, which is the probability of $X$ being at most $x$; thus $$F(x) \ := \ \int_{-\infty}^x f(t)dt.$$

Consider the cumulative distribution function of a normal with mean zero and variance 1/2. Thus $$f(x) \ = \ \frac1{\sqrt{\pi}} e^{-x^2}, \ \ \ F(x) \ = \ \frac1{\sqrt{\pi}} \int_{-\infty}^x e^{-t^2} dt;$$ the variance is taken to be 1/2 to make the argument of the exponential $-t^2$ and not $-t^2/2$. As we don't have a closed form anti-derivative of $e^{-t^2}$, we're in trouble. We could try to use the series expansion for $e^u$ and interchange the integral and the sum, but as the lower bound is $-\infty$ the resulting expression contains alternating signed infinities. Fortunately we know $F(0) = 1/2$ as our density is symmetric about zero.

So if $x \ge 0$ (the case $x < 0$ is handled similarly with a subtraction now) we can write $$F(x) \ = \ \frac12 + \frac12\erf(x),$$ where (in a lovely set of normalizations!) $$\erf(x)\ := \ \frac{2}{\sqrt{\pi}}\int_0^x e^{-t^2}dt.$$ While it would've been disastrous to use $e^u = \sum_n u^n/n!$ when integrating from $-\infty$, now we can insert it with $u = -t^2$ and interchange the integral and the sum (see for example \cite{La} for a justification), yielding
\begin{eqnarray}
   \erf(x) &\ := \ &  \frac{2}{\sqrt{\pi}}\int_0^x e^{-t^2} dt \nonumber\\ &=& \frac{2}{\sqrt{\pi}}\int_0^x \sum_{n=0}^\infty \frac{(-t^2)^n}{n!} dt \nonumber\\  &=& \frac{2}{\sqrt{\pi}} \sum_{n=0}^\infty \frac{(-1)^n}{n!} \int_0^x t^{2n} dt \nonumber\\ &=&   \frac{2}{\sqrt{\pi}}\sum_{n=0}^\infty \frac{(-1)^n}{n!} \frac{t^{2n+1}}{2n+1}\Bigg|_0^x \nonumber\\ &=&  \frac{2}{\sqrt{\pi}}\sum_{n=0}^\infty \frac{(-1)^n x^{2n+1}}{n! (2n+1)}, \nonumber
\end{eqnarray} which converges for all $x$!

The point of the above is that we have to really think what it means to solve a problem. We have a very rapidly converging series that we've named, but it is an infinite sum. If this bothers you, think back to trig functions such as sine, cosine and tangent -- we consider those nice, but to find a general value for any of these also requires evaluating an infinite sum.


\ \\

\end{document}